\documentclass[12pt,a4paper]{amsart} \usepackage{amssymb,amscd} \usepackage{xypic}
\usepackage[english]{babel}
\theoremstyle{plain}
\newtheorem{theorem}{Theorem}[section]
\newtheorem{lemma}[theorem]{Lemma}
\newtheorem{proposition}[theorem]{Proposition}
\newtheorem{corollary}[theorem]{Corollary}
\theoremstyle{definition}
\newtheorem{definition}[theorem]{Definition}
\newtheorem{remark}[theorem]{Remark}
\newtheorem{remarks}[theorem]{Remarks}
\newtheorem{example}[theorem]{Example}

\numberwithin{equation}{section}

\newcommand\cA{{\mathcal A}}

\newcommand\cO{{\mathcal O}}

\newcommand\kb{\bar{k}}

\newcommand\aff{\operatorname{aff}}
\newcommand\qa{\operatorname{qa}}

\newcommand\End{\operatorname{End}}
\newcommand\Adit{\operatorname{Addit}}
\newcommand\Aut{\operatorname{Aut}}
\newcommand\Hom{\operatorname{Hom}}

\newcommand\Pic{\operatorname{Pic}}
\newcommand\Spec{\operatorname{Spec}}

\newcommand\Imm{\operatorname{Imm}}

\def\res #1{\setbox1=\hbox{\text{$\left .\vbox
 to9pt{}\right\vert_{#1}$}} \lower 2pt\hbox to
 \wd1{\text{$\left .\vbox to9pt{}\right\vert_{#1}$}}}

\def\punto{\cdot}

\title[Principal bundles and quasi-abelian varieties]{Principal bundles, quasi-abelian varieties and structure
of algebraic
groups}

\author{Carlos Sancho de Salas}
\address{Universidad de Salamanca\\ Departamento de Matem\'aticas\\
Plza. de la Merced s/n, Salamanca\\ 37008, Espa\~na}
\email{mplu@usal.es}
\author{Fernando Sancho de Salas}
\address{Universidad de Salamanca\\ Departamento de Matem\'aticas\\
Plza. de la Merced s/n, Salamanca\\ 37008, Espa\~na}
\email{fsancho@usal.es}

\begin{document}

\begin{abstract} We classify principal bundles over anti-affine
schemes with affine and commutative structural group. We show that
this yields the classification of quasi-abelian varieties over a
field $k$ (i.e., group $k$-schemes $G$ such that $\mathcal
O_G(G)=k$). The interest of this result is given by the fact that
the classification of smooth group $k$-schemes is reduced to the
classification of quasi-abelian varieties and of certain affine
group schemes.

\end{abstract}

\maketitle

\tableofcontents

\section{Introduction}
\label{sec:introduction}

Let $k$ be a field and $G$ a group $k$-scheme  of finite type. We
say that $G$  is a quasi-abelian variety if $\cO_G(G)=k$. Examples
include abelian varieties, their universal vector extensions (in
characteristic $0$ only) and certain semi-abelian varieties. The
main motivation to study quasi-abelian varieties is the fact that
the classification of group schemes over fields is essentially
reduced to the classification of quasi-abelian varieties and of
affine group schemes. In fact one has (Theorem \ref{genClass}):

\smallskip

\begin{theorem} {\bf (Structure of algebraic
groups)} Every connected smooth  $k$-scheme in groups $G$
decomposes as
$$G\simeq (\overline
G\times\mathcal A)/H$$ where $\overline G$ is an affine connected
 group without finite quotients, $\mathcal A$ is a quasi-abelian
variety and $H$ is an affine commutative group $k$-scheme
satisfying:

- $H$ is contained in the center of $\overline G$.

- $\mathcal A_{\aff}\subset H\subset \mathcal A$ and $H/\mathcal
A_{\aff}$ is finite, with $\mathcal A_{\aff}=$  affine part of
$\cA$.

- $H$ is submerged in $\overline G\times\mathcal A$ through the
diagonal morphism.

This decomposition is unique up to isomorphisms of $\overline G$
and $\mathcal A$.
\end{theorem}

This theorem is essentially contained in the work of Rosenlicht
\cite{Ro56} over an algebraically closed field. One can extend it
to arbitrary  fields using the results of \cite{BLR90}. We include
a proof in order to be self-contained.

%

This result reduces the classification of algebraic groups to the
classification of  affine groups and quasi-abelian varieties and
motivates the aim of this paper: the structure and classification
of quasi-abelian varieties. A second motivation comes from the
problem of classification of homogeneous varieties. This problem
is essentially fulfilled in the proper case (see \cite{Sa03}). The
next step is to deal with the anti-affine case (anti-affine means
that the variety has only constant global functions). This case
seems accessible because these varieties are rigid (as we show in
Theorem \ref{lemmarigidity}). It is convenient to study first the
case of groups. Firstly, because they are a particular case of
homogenous variety. Secondly, because this study should be useful
to understand the structure of the automorphism group of these
varieties (notice that in the proper case this group is almost
classifying).

Despite its interest, the study of quasi-abelian varieties is
limited in the literature; they only appear implicitly in work of
Rosenlicht and Serre (see \cite{Ro58, Ro61, Se58a}). In Analytic
Geometry there exists a notion of quasi-abelian variety (see
\cite{AK01}) which is stronger than the algebraic one. This means
an algebraic variety over  $\Bbb C$ such that, as an analytic
variety, has no non constant global functions. Clearly these
varieties are quasi-abelian in the algebraic sense, but the
converse is not true. For example, the universal vectorial
extensions of abelian varieties are quasi-abelian in the algebraic
sense but they have non constant analytical global functions
because they are Stein.

Here we obtain the structure of quasi abelian varieties and we
reduce their classification to that of abelian varieties.

With respect to the structure of quasi-abelian varieties one first
notices that Chevalley's theorem implies that a quasi-abelian
variety is a principal bundle over an abelian variety $A$ with
affine, commutative and connected structural group $G$. We shall
prove that the  classification of quasi-abelian varieties as
groups is equivalent to their classification as principal bundles.
That is, two quasi-abelian varieties are isomorphic (as group
schemes) if and only if they are isomorphic as principal bundles
over isomorphic abelian varieties with isomorphic structural
groups (see Theorem \ref{fibradoGrupo} and Corollary
\ref{prin-quasi}). This will be a consequence of the rigidity of
quasi-abelian varieties. In this direction, we shall give a
general rigidity theorem for anti-affine schemes, which, as
mentioned above, has its interest in the classification of
anti-affine homogeneous varieties.

Next we deal with the classification of principal bundles over an
anti-affine scheme $Y$ with affine and commutative structural
group $G$. We shall always assume that a principal bundle has a
rational point (see Remark \ref{rational}). Let us denote
$\operatorname{Prin}(G,Y)$ the set of isomorphism classes of
principal $G$-bundles over $Y$ and
$\operatorname{Prin}(G,Y)_{ant}$ the set of isomorphism classes of
anti-affine principal $G$-bundles over $Y$. If $Y$ is an abelian
variety, let us denote $\operatorname{Prin}(G,Y)_{ant}^{st}$ the
set of isomorphism classes of anti-affine principal $G$-bundles
over $Y$ which are stable under translations on $Y$ (see
Definition \ref{trans-inv}). Theorem \ref{fibradoGrupo} says that
the quotient of $\operatorname{Prin}(G,Y)_{ant}^{st}$ by the
automorphism group of $G\times Y$ coincides with the set of
isomorphism classes of quasi-abelian varieties with affine part
isomorphic to $G$ and abelian part isomorphic to $Y$.

The key point for our classification of principal bundles will be
its relation with the Cartier dual of $G$ and the Picard scheme of
$Y$, that we explain now. Let $\pi\colon P\to Y$ be a principal
$G$-bundle. Each character $\chi$ of $G$ determines an invertible
subsheaf $\mathcal L_{\chi}$ of $\pi_*\mathcal O_P$, namely the
subsheaf of functions of $P$ over which $G$ acts by that
character; hence, the principal $G$-bundle $\pi\colon P\to Y$
defines a morphism of functors of groups $G^D\to \bold{Pic} (Y)$,
where $G^D$ is the Cartier dual (functor) of $G$. We shall prove
that this morphism classifies the bundle (see Theorem
\ref{clasFib} for the precise statement). Once
$\operatorname{Prin}(G,Y)$ is determined, we deal with
$\operatorname{Prin}(G,Y)_{ant}$ and
$\operatorname{Prin}(G,Y)_{ant}^{st}$ (see Theorems
\ref{antiafin-criterio}, \ref{aaImm} and \ref{clasification-1}).

From here, making use of the knowledge of $G^D$ for either a
unipotent or a multiplicative type $G$ and the structure of
$\bold{Pic} (Y)$, we shall obtain a full description of
$\operatorname{Prin}(G,Y)$, $\operatorname{Prin}(G,Y)_{ant}$ and
$\operatorname{Prin}(G,Y)_{ant}^{st}$  (see Theorems
\ref{multiplicative}, \ref{vectorial}, \ref{vectExt} and
\ref{general}). In particular, we obtain  the known classification
theorems of principal bundles over an abelian variety whose
structural group is either a vector space or the multiplicative
group (see \cite{MM74, Se59, Ro58}). This ``Cartier-perspective''
will be also very useful for the classification of anti-affine
homogenous varieties, since it is not difficult to prove that
these varieties are principal bundles over proper homogeneous
varieties.

From this perspective we obtain our main result (Theorem
\ref{classqa}) that classifies quasi-abelian varieties
 over an arbitrary field $k$:

\begin{theorem} Let us denote $k_s$ the separable closure of $k$. Then {\sl
to give a quasi-abelian variety $\cA$ over $k$ with affine part
$G$ and abelian part $Y$ is equivalent to give the following data:
\begin{enumerate}
\item A sublattice $\Lambda \subset \Pic^0(Y_{k_s})$, stable under the
action of the Galois group $\mathcal G({k_s}/k)$. \item A linear
subspace $V\subset H^1(Y,\mathcal O_Y)$,
\end{enumerate} such that $\Lambda\simeq X(G_{k_s})$ and $V\simeq
\Adit(G)$, where $\Adit(G)$ is the vector space of additive
functions of $G$ and $X(G_{k_s})$ is the group of characters of
$G_{k_s}$. These data are given up to group automorphisms of $Y$}.
\end{theorem}

This  classification  was obtained in \cite{Sa01}, with similar
techniques, when $k$ is an algebraically closed field. It has also
been  proved independently by M. Brion (see \cite[Thm.~2.7]{Br}).

\smallskip

As a consequence of the classification theorem we obtain   that
every quasi-abelian variety over a field of positive
characteristic is semi-abelian. One also obtains that, over an
arbitrary base field, the affine part of a quasi-abelian variety
is smooth.

\bigskip

\noindent {\bf Acknowledgements.} We wish to thank M. Brion for
his patient and enlightening attention to this paper. His valuable
remarks have allowed this paper to reach its final form.
\bigskip

\noindent {\bf Notation and conventions.} Throughout this article,
$k$ is a field with separable closure $k_s$ and algebraic closure
$\kb$.

By a \emph{scheme}, we mean a scheme of finite type over $k$,
unless otherwise specified; a point of a scheme will always mean a
valued point. Morphisms of schemes are understood to be
$k$-morphisms, and products are taken over $k$. A \emph{variety}
is a separated and geometrically integral scheme. A {\sl functor}
is always a functor from the category of $k$-schemes (or
$k$-algebras) to the category of sets. The functor of points of a
scheme $X$ is still denoted by $X$.

As in \cite{Br} we say that a scheme $X$ is anti-affine if
$\mathcal O_X(X)=k$.

We shall use a boldface type to denote functors like $\bold{Aut}$,
$\bold{Pic}$, $\bold{Hom}$, etc (functor of automorphisms, Picard
functor, functor of homomorphisms, etc) and for the schemes
representing them (when they exist). We shall use a non boldface
type like $\Aut$, $\Pic$, $\Hom$, etc for the sets of
automorphisms, Picard group, homomorphisms, etc.

 By an \emph{algebraic group}  we mean a {\it smooth}
group scheme $G$, possibly non-connected. An \emph{abelian
variety} is a connected and complete algebraic group. For these,
we refer to \cite{Mu70}, and to \cite{Bo91} for affine algebraic
groups. For any group scheme $G$, a $G$-scheme means a scheme
endowed with an action of $G$ on it. A group $G$ is   {\sl of
multiplicative type} if $G_{\overline k}$ is diagonalizable. A
{\sl torus} is a  smooth group of multiplicative type.

For any group $G$, $X(G)$ denotes the group of characters of $G$,
i.e., $X(G)=\Hom_{groups}(G,G_m)$.

It is well known that any commutative affine group $G$ has a
unique multiplicative type subgroup $\mathcal K$ such that
$G/\mathcal K= U$ is unipotent. We say that $\mathcal K$ (resp. $
U$) is the {\sl multiplicative type part  of $G$} (resp. the {\sl
unipotent part of} $G$). It is not true in general that   $G=
U\times\mathcal K$, but it holds when $k$ is perfect.

For any connected group scheme $G$  we denote by $G_{\aff}$ the
smallest normal connected affine subgroup such that the quotient
$G/G_{\aff}$ is an abelian variety. We shall call $G_{\aff}$
(resp. $G/G_{\aff}$)  the {\sl affine part of }$G$ (resp. the {\sl
abelian part of} $G$). The existence of $G_{\aff}$ is due to
Chevalley in the setting of algebraic groups over algebraically
closed fields; in this case $G_{\aff}$ is an algebraic group as
well, see \cite{Ro56,Ch60}. Chevalley's theorem easily implies the
existence of $G_{\aff}$ for any connected group scheme $G$, see
\cite[Lem.~IX.2.7]{Ra70} or \cite[Thm.~9.2.1]{BLR90}. If $G$ is an
algebraic group and $k$ is perfect, then $G_{\aff}$ is also an
algebraic group. If $k$ is not perfect, then  $G_{\aff}$ is
connected but it might be non smooth. We do not know if $G_{\aff}$
can be non-reduced. In any case, it is immediate that $G_{\aff}$
is quasi-reduced. By this we mean

\begin{definition}\label{qreduced}
We say that a  group scheme $G$ is quasi-reduced if for any
subgroup  $H\subset G$  such that $H_{red}=G_{red}$ one has $H=G$.
If $G$ is connected, this is equivalent to say  that $G$ does not
admit finite quotients.
\end{definition}

\begin{remark}\label{ToroQreduc} Let $G$ be a group of multiplicative type. Then, for any
$n\in\mathbb N$, the multiplication  $G\overset{\cdot n}\to G$ is
an isogeny. Moreover, if $n=\vert G_{\overline k}/(G_{\overline
k})_{\text{red}}\vert $, then $n G$ is smooth and connected. Hence
$nG$ coincides with the reduced and connected component at the
origin of $G$. In conclusion, {\sl if $G$ is a connected and
quasi-reduced group of multiplicative type, then it is a torus}.
\end{remark}

\section{Quasi-abelian part of a group scheme. Basic properties of
quasi-abelian varieties: Rigidity}

In this section we establish known results about  quasi-abelian
varieties and we generalize the rigidity theorem of proper
varieties to anti-affine schemes.

The following results, stated without proof, can be found in
\cite[Sec.~III.3.8]{DG70}.

\begin{theorem} If $G$ is a quasi-abelian variety then it is smooth  and connected.
\end{theorem}

If $G$ is a group scheme, then $A=H^0(G,\mathcal O_G)$ is a Hopf
$k$-algebra and one has a natural morphism of groups:
$$\pi_{\text{aff}}\colon G\to \operatorname{Aff}(G)$$ where
$\operatorname{Aff}(G)=\Spec A$.

This affine group $\operatorname{Aff}(G)$ is called the {\sl
affinization group of $G$} and it satisfies trivially the
universal property:
$$\Hom_{groups}(G,H)=\Hom_{groups}(\operatorname{Aff}(G),H)$$ for
any affine group $H$.

\begin{definition} For each group scheme $G$ we denote
$G_{\qa}=\ker\pi_{\text{aff}}$ and we call it {\sl the
quasi-abelian part} of $G$. One has
$G/G_{\qa}=\operatorname{Aff}(G)$.
\end{definition}

\begin{proposition} The quasi-abelian part of $G$ is a
quasi-abelian variety.
\end{proposition}

\begin{theorem} Let $G$ be a quasi-abelian variety and $H$ a
connected  group. If $f\colon G\to H$ is a morphism of schemes
such that $f(e)=e$, then
\begin{enumerate}
\item $f$ is a morphism of groups,
\item $f$ takes values in the center of $H$,
\item $f$ takes values in $H_{\qa}$.
\end{enumerate}
\end{theorem}

\begin{theorem}\label{structqa} If $G$ is a quasi-abelian variety
then its group structure is unique (once the neutral point is
fixed) and it is commutative. Moreover if $G$ is a subgroup of a
group $H$, then  it is contained in the center of $H$.
\end{theorem}

The latter two theorems can be easily obtained from the rigidity
theorem for anti-affine schemes that we shall next prove. It
generalizes the rigidity theorem of abelian varieties and it shows
that rigidity is not as much a consequence of properness but of
anti-affinity.

\begin{lemma} Let $X$ be an anti-affine scheme and $Y$ an affine scheme. Any morphism of schemes
$X\to Y$ is constant (i.e., it factors through a morphism $\Spec
k\to Y$).
\end{lemma}

\begin{proof} Obvious.
\end{proof}

\begin{theorem}{{\bf (Rigidity of
anti-affine schemes)}}\label{lemmarigidity} Let $X$, $Y$ and $Z$
be schemes, $X$ anti-affine with some rational point, $Y$
connected and $Z$ separated. Let
  $$ f\colon X\times Y \longrightarrow Z
$$ be a morphism. If there exist a closed point $y_0\in Y$ such
that $f_{|X\times \{y_0\}}$ is a constant morphism, then $f$
factors
  $$
\xymatrix{ X\times Y \ar[r]^{f} \ar[d]_{p_2} & Z \\ Y \ar[ur]_{g}&
}
$$
  where $p_2$ is the second projection.
\end{theorem}

\begin{proof} We shall fix a rational point $x_0\in X$. Let us
define $g\colon Y\to Z$ as $g(y)=f(x_0,y)$. We claim that
$f=g\circ p_2$.

(a) Assume that $Z$ is an affine scheme, $Z=\Spec A$. Then $f$ is
constant on $X$, because to give a morphism $X\times Y\to Z$ is
equivalent to give a morphism of $k$-algebras $A\to H^0(X\times
Y,\mathcal O_{X\times Y})= H^0(Y,\mathcal O_Y)$, i.e., a morphism
$Y\to Z$.

(b) If the morphism $f_0\colon (X\times Y)_{top}\to Z_{top}$,
between the underlying topological spaces, factors through
$g_0\colon Y_{top}\to Z_{top}$ (i.e. $f_0=g_0\circ (p_2)_0$), then
$f$ factors. Indeed, for each affine open subscheme $U\subset Z$,
let $V=g_0^{-1}(U)$. One has $f_0^{-1}(U)=X\times V$. Then $f$
maps $X\times V$ into $U$ and the morphism $f\colon X\times V\to
U$ factors through $g\colon V\to U$ (by (a)). So if
$Z=\bigcup_iU_i$ is an affine open covering, then $X\times
Y=\bigcup_i f^{-1}(U_i)$ is an open covering and $f$ factors over
each $f^{-1}(U_i)$.

(c) We can assume that $Y$ is irreducible. Indeed, let
$Y=Y_0\cup\cdots\cup Y_n$ be a decomposition on irreducible
components such that $y_0\in Y_0$. Let $Y_i$ be another component
 meeting $Y_0$. If the claim holds when $Y$ is an irreducible
scheme, then $f$ is constant along fibers over $Y_0$. So, $f$ is
constant along fibers over $Y_0\cap Y_i$, and then along fibers
over $Y_0\cup Y_i$. By recurrence, $f$ is constant along fibers
over the whole $Y$.

Now let $T\subset X\times Y$ be the sub-scheme of points $t$ such
that $f(t) =(g\circ p_2)(t)$. Since $Z$ is separated, $T$ is a
closed sub-scheme.

(d) $T$ contains a open neighborhood of $X\times\{ y_0\}$. Indeed,
let $\mathcal O$ be the local ring of $Y$ at $y_0$, $\mathfrak m$
its maximal ideal and let us denote $X_n=X\times\Spec \mathcal
O/\mathfrak m^n\subset X\times Y$. It is clear that $f(X_n)$ is a
finite subscheme of $Z$ (supported on $z_0$). Then $f(X_n)$ is an
affine scheme and, by (a), $f\res {X_n}$ factors through $\Spec
\mathcal O/\mathfrak m^n$, i.e. it is equal to $g\circ p_2$. Hence
$T\supset X_n$ for all $n$. Since $\bigcap_n \mathfrak m^n=0$, we
conclude that $T$ contains a neighborhood of $X\times\{ y_0\}$ in
$X\times Y$.

Now, since $Y$ is irreducible, each irreducible component of
$X\times Y$ maps surjectively on $Y$. So, all of them cut
$X\times\{ y_0\}$. By (d) $T$ contains a non empty open subset of
each one. Since $T$ is closed, it contains all irreducible
components of $X\times Y$. So $T_{top}=(X\times Y)_{top}$ and we
conclude by (b).
\end{proof}

\section{Structure of algebraic
groups}

We  give a structure theorem for algebraic groups that sums up
results of Chevalley, Rosenlicht, Demazure-Gabriel and
\cite{BLR90}.

\begin{theorem}\label{genClass} {\bf (Structure of algebraic
groups)} Every connected algebraic group $G$ decomposes as
$$G\simeq (\overline
G\times\mathcal A)/H$$ where $\overline G$ is an affine connected
quasi-reduced group (see definiton \ref{qreduced}), $\mathcal A$
is a quasi-abelian variety and $H$ is an affine commutative group
scheme satisfying:

- $H\subset Z(\overline G)$.

- $\mathcal A_{\aff}\subset H\subset \mathcal A$ and $H/A_{\aff}$
is finite.

- $H$ is submerged in $\overline G\times\mathcal A$ through the
diagonal morphism.

This decomposition is unique up to isomorphisms of $\overline G$
and $\mathcal A$.
\end{theorem}

\begin{proof}
If we denote $\overline G=G_{\aff} , \mathcal A=G_{\qa},
H=G_{\aff}\cap G_{\qa}$, then one has the desired decomposition.
Indeed: the quotient of $G$ by $G_{\aff}\cdot G_{\qa}$ is trivial
because it is a quotient of the abelian variety $G/G_{\aff}$ and a
group quotient of the affine group $G/G_{\qa}$ and so it is an
abelian variety and an affine group. Hence $G=G_{\aff}\cdot
G_{\qa}$. Moreover $\mathcal A/H\hookrightarrow G/G_{\aff}$ is
abelian and so $\mathcal A_{\aff}\subset H$
 and $H/\mathcal A_{\aff}\subset \mathcal A/\mathcal A_{\aff}$ is
closed and affine (because $H\subset G_{\aff}$ is affine) and then
it is finite.

 Conversely, if
$G\simeq (\overline G\times\mathcal A)/H$ as in the theorem
hypothesis, then $\overline G$ and $\mathcal A$ are normal
connected subgroups of $G$, $H=\overline G\cap \mathcal A$,
$\overline G$ is affine quasi-reduced and $\mathcal A$ is a
quasi-abelian variety. Moreover $G/\overline G\simeq \mathcal A/H$
is an abelian variety (because $\mathcal A/H$ is a quotient of
$\mathcal A/\mathcal A_{\aff}$, an abelian variety) and
$G/\mathcal A$ is affine because it is a quotient of $\overline
G$. Hence $\overline G= G_{\aff},\mathcal A= G_{\qa}$ and then
$H=G_{\aff}\cap G_{\qa}$.
\end{proof}

This theorem says that the classification of algebraic groups is
essentially reduced to the classification of affine groups and
quasi-abelian varieties.

We can refine this result when the base field is perfect  in the
following way (see also \cite{Br}, sections 3.2 and 3.3, for
related results):

\begin{proposition}\label{prop:preari}
Let $G$ be a connected algebraic group over a perfect field $k$.
Then there exists a reduced, connected and affine group
$\widetilde G$, a quasi-abelian variety $\cA$ and an isogeny
$$\phi\colon (\widetilde G\times \mathcal A)/\mathcal U\to G$$ such
that $\phi\res{\widetilde G}$ and $\phi\res{\mathcal A}$ are
injective morphisms, where $\mathcal U$ is the unipotent part of
$\mathcal A_{\aff}$
 and $\mathcal U\to \widetilde G\times \mathcal A$ is the
diagonal morphism induced by an immersion $\mathcal U
\hookrightarrow Z(\widetilde G)$. Moreover, with these conditions,
$\widetilde G$ and $ \mathcal A$ are unique up to isomorphisms. In
fact $\mathcal A\simeq G_{qa}$ and $\widetilde G$ is a
quasi-complement of the multiplicative part of $ \mathcal
A_{\aff}$ in $G_{\aff}$.
\end{proposition}

\begin{proof}
Let us take $\mathcal A=G_{qa}$ and let us denote by $S$ the
multiplicative part of $ \mathcal A_{\aff}$. By Theorem
\ref{genClass} it suffices to show that $S$ has a quasi-complement
in $G_{\aff}$. This is well known if $G_{\aff}$ is reductive. For
the general case, let $G'$ be a quasi-complement of $S$ in
$G_{\aff}/R_u$, where $R_u$ is the unipotent radical of $
G_{\aff}$. If $\pi\colon G_{\aff}\to G_{\aff}/R_u$ is the quotient
map, then $\widetilde G=\pi^{-1}(G')$ is a quasi-complement of $S$
in $ G_{\aff}$.

The uniqueness of $\widetilde G$ and $ \mathcal A$ is not
difficult.\end{proof}

\section{Quasi-abelian varieties as principal bundles}

 As we have seen, a quasi-abelian variety $\cA$ is a commutative
group (Theorem \ref{structqa}). Moreover there exists a connected
and affine subgroup $G\subset \cA$ such that the quotient $\cA/G$
exists and it is an abelian variety (Chevalley's structure
theorem). Therefore a quasi-abelian variety may be thought of as
an extension of an abelian variety by an affine commutative group,
or as a principal bundle on an abelian variety with affine and
commutative structural group. Recall that a principal bundle over
a scheme $Y$ with structural group $G$  is a $G$-scheme $P$
together with a morphism of $G$-schemes $P\to Y$ (where $G$ acts
trivially on $Y$) such that the natural map
\[\aligned  G\times P&\to P\times_Y P\\ (g,p)&\mapsto (g\cdot
p,p)\endaligned\] is an isomorphism. For short, we say that $P\to
Y$ is a principal $G$-bundle.

\begin{remark}[Extra hypothesis] \label{rational} We shall always assume that a principal $G$-bundle
$P$ over $Y$ has a rational point, since this is the case when $P$
is a quasi-abelian variety. As we shall see, this  implies (in our
hypothesis, i.e., $G$ a commutative affine group and $Y$ an
anti-affine scheme with some rational point) that a principal
$G$-bundle over $Y$ is locally split: there exists a Zariski open
covering $U_i$ of $Y$ such that $P_{\vert U_i}=U_i\times G$. This
is why we have used the terminology of {\it principal bundles}
(which is more common in differential geometry) instead of {\it
torsors}.
\end{remark}

A morphism $f\colon P\to P'$ of principal $G$-bundles over $Y$ is
a morphism of $G$-schemes over $Y$.

We denote by $\operatorname{Prin}(G,Y)$ the set of isomorphism
classes of principal $G$-bundles over $Y$ and by
$\operatorname{Prin}(G,Y)_{ant}$ the set of isomorphism classes of
anti-affine principal $G$-bundles over $Y$. If $Y$ is an abelian
variety, we shall denote by $\operatorname{Prin}(G,Y)_{ant}^{st}$
the set of isomorphism classes of anti-affine principal
$G$-bundles over $Y$ which are stable under translations on $Y$
(see Definition \ref{trans-inv}).

It is clear that $\Aut_{groups}(G)$ and $\Aut_{schemes}(Y)$ act on
$\operatorname{Prin}(G,Y)$,  $\operatorname{Prin}(G,Y)_{ant}$ and
$\operatorname{Prin}(G,Y)_{ant}^{st}$.

We say that two quasi-abelian varieties are isomorphic if they are
isomorphic as group schemes. Two isomorphic quasi-abelian
varieties have isomorphic affine parts and isomorphic abelian
parts. We shall denote by $\operatorname{Quasiabel}(G,Y)$ the set
of isomorphism classes of quasi-abelian varieties whose affine
part is isomorphic to $G$ and whose abelian part is isomorphic to
$Y$. The aim of this section is to prove that
\[ \operatorname{Prin}(G,Y)_{ant}^{st}/\Aut_{groups}(G\times Y) =
\operatorname{Quasiabel}(G,Y)\] The key point is to show that if
$P$ is an anti-affine principal $G$-bundle over an abelian variety
$Y$ and it is stable under translations on $Y$, then $P$ admits a
(essentially unique) group structure such that $P$ is a
quasi-abelian variety with affine part $G$ and abelian part $Y$.
This will be done in Theorem \ref{fibradoGrupo}.

\begin{lemma}\label{rigAut}
Let $G$ be a commutative group scheme and $\pi\colon P\to Y$  a
principal $G$-bundle. Let us denote $\bold{Aut}_Y^G(P)$ the
functor of automorphisms of principal $G$-bundles of $P$. One has
$$\bold{Aut}_Y^G(P)=\bold{Hom}_{schemes}(Y,G)$$
In particular, if $G$ is affine and $Y$ is anti-affine, then
$\bold{Aut}_Y^G(P)=G$.
\end{lemma}

\begin{proof}
Since $G$ is commutative, it is clear that
$\bold{Aut}_{G-schemes}(G)=G$ and then
$\bold{Aut}_{G-schemes}(Z)=G$ for every $G$-scheme $Z$ on which
$G$ acts free and transitively. Then one has a morphism
$$\aligned \bold{Aut}_Y^G(P)&\to\bold{Hom}_{k-schemes}(Y,G)\\ \tau &\mapsto
f_{\tau}\endaligned$$ where $f_{\tau}(y)$ is the automorphism of
$G$ induced by $\tau$ in the fiber of the (valued) point $y$.
Conversely, given $f\colon Y\to G$, one has a $G$-automorphism
$\tau_f\colon P\to P$, $\tau_f(p)=f(p)\cdot p$. We conclude
immediately.
\end{proof}

\begin{definition}\label{trans-inv} Let $Y$ be a group scheme and $G$ an
affine commutative group. A principal $G$-bundle $\pi\colon P\to
Y$ is said to be stable under translations on $Y$ if for each
point $y\colon Z\to Y$ there exists a faithfully flat base change
$Z'\to Z$ and a morphism of $G$-schemes $\varphi_y \colon P\times
{Z'}\to P\times {Z'}$ such that the diagram:
$$\CD P\times {Z'} @>{\varphi_y}>> & P\times {Z'}
  \\ @V{\pi }VV & @V{\pi }VV
\\Y\times {Z'} @>{\tau_y}>> & Y\times {Z'}\endCD $$ is commutative,
where $\tau_y$ is the translation by $y$.\end{definition} More
briefly, a principal $G$-bundle $P\to Y$ is stable under
translations on $Y$ if any translation on $Y$ extends (up to a
faithfully flat base change) to an automorphism of $G$-schemes of
$P$.

For example, if $\cA$ is a quasi-abelian variety with affine part
$G$ and abelian part $Y$, then $\cA$ is a principal $G$-bundle
over $Y$ and it is obviously stable under translations on $Y$. We
now see that the converse also holds.

\begin{theorem}\label{fibradoGrupo} Let $Y$ be an abelian variety,
$G$ an affine commutative group scheme and $\pi\colon P\to Y$ a
principal $G$-bundle. Then $P\to Y$ is stable under translations
on $Y$ if and only if $P$ admits a group structure such that:
\begin{enumerate}
\item[(i)] $\pi\colon P\to Y$ is a morphism of groups,

\item[(ii)] the kernel of $\pi$ is isomorphic to  $G$ as a $G$-scheme,
and

\item[(iii)] the translations by points of $P$ commute with the action of
$G$.
\end{enumerate}

Moreover, this group structure is unique (once the neutral point
on the fiber of $0\in Y$ is fixed), and it is commutative. If in
addition $P$ is anti-affine, then it is a quasi-abelian variety.
\end{theorem}

\begin{proof} Assume that $P$ has a group structure satisfying (i)-(iii). First notice that $P$ is
commutative; indeed, let $G_0,P_0$ be the connected components
through the origin of $G,P$, respectively. It is clear that
$G\cdot P_0=P$ and then it is enough to prove that $P_0$ is
commutative. So, replacing $P, G$ by $P_0,G_0$, we can suppose
that $P$ is connected. On the one hand the quotient of $P$ by its
quasi-abelian part is affine and then the quotient by its center
subgroup is also affine; on the other hand this quotient is a
quotient of $P/G=Y$ (because $G$ is in the center of $P$) and then
it is proper. Hence the quotient of $P$ by its center is trivial
and $P$ is commutative. Now let us see that $\pi\colon P\to Y$ is
stable under translations on $Y$, i.e., each translation on $Y$
lifts to an automorphism of $G$-schemes on $P$ (after a faithfully
flat base change). Indeed, since $P\to Y$ is a faithfully flat
morphism, each point $y$ of $Y$ has some point in its fiber by
$\pi$ (after a faithfully flat base change). So it is enough to
define on $P$ the translation morphism by any point of this fibre.

Assume now that $P$ is stable under translations on $Y$.  Let
$\mathbf{Aut}^Y(P/Y)$ be the functor $\mathbf{Aut}^Y(P/Y)(Z)
=\{$automorphisms $\varphi\colon P_Z\to P_Z$ of $G$-schemes which
descend to a translation on $Y_Z\}$. One has an exact sequence of
functors of groups:
$$0\to G\to \mathbf{Aut}^Y(P/Y)\overset p\to Y\longrightarrow 0$$
where $p$ is the morphism that maps each automorphism $\varphi$ to
the induced translation on $Y$. The surjectivity of $p$ (for the
faithfully flat topology) is due to the hypothesis, i.e.,
$\pi\colon P\to Y$ being stable under translations, and the kernel
of $p$ is $G$ by Lemma \ref{rigAut}. $\mathbf{Aut}^Y(P/Y)$ acts
freely on $P$. Moreover this action is transitive: indeed, given
two points $p_1,p_2$ of $P$ there exists a translation on $Y$
transforming $\pi (p_1)$ on $\pi (p_2)$, so we can assume that
$\pi (p_1)=\pi (p_2)$. One concludes the transitivity because $G$
acts transitively on the fibres of $\pi$. Now let us fix a
rational point $e\in \pi^{-1}(0)$. Transforming $e$ by
$\mathbf{Aut}^Y(P/Y)$ we obtain that $\mathbf{Aut}^Y(P/Y)\simeq P$
and so $P$ has a group structure satisfying the required
conditions.

{\it Uniqueness:} the translations on $P$ define a group immersion
$P\hookrightarrow \mathbf{Aut}^Y(P/Y)$, whose composition with the
isomorphism $\mathbf{Aut}^Y(P/Y)\simeq P$ is the identity. So the
group structure of $P$ is the one induced by the isomorphism
$\mathbf{Aut}^Y(P/Y)\simeq P$.\end{proof}

\begin{corollary}\label{prin-quasi} Two quasi-abelian varieties are isomorphic (as
groups) if and only if their affine parts and their abelian parts
are respectively isomorphic and they are isomorphic as principal
bundles. In other words, one has a bijection
\[ \operatorname{Prin}(G,Y)_{ant}^{st}/\Aut_{groups}(G\times Y) =
\operatorname{Quasiabel}(G,Y)\]
\end{corollary}


\begin{remark} As we have seen in the proof of  Theorem
\ref{fibradoGrupo}, the existence and the uniqueness of the group
structure of a principal $G$-bundle over a group $Y$ only needs
that $\Hom_{schemes}(Y,G)=G$; that is, it only needs that any
morphism of schemes $Y\to G$ is constant. Hence Theorem
\ref{fibradoGrupo} can be extended to different cases. For
example, for the calculation of the extensions of unipotent groups
(smooth and connected but possibly non commutative) by
multiplicative type groups. In particular, this would reduce the
classification of affine abelian groups (over an arbitrary field)
to the classification of unipotent groups and of their principal
bundles with multiplicative type structural group.

\end{remark}

\section{Cartier Dual and classification
of principal bundles}

In this section we obtain the classification of principal
$G$-bundles over an anti-affine scheme $Y$, with $G$ an affine
commutative group scheme. It generalizes well known results about
the subject in the particular cases when the structural group $G$
is either a torus or a vector space (see \cite{MM74, Se59, Ro58}).
Moreover this result allows us to see that the differences between
these cases (torus and vector space) come only from the different
structure of the respective Cartier dual groups (local and
discrete, respectively) \

\subsection{$i$-component of linear representations}

Let $G=\Spec A$ be an affine group $k$-scheme. Let us denote
\[ I =\text{ set of finite sub-coalgebras of } A.\] For each $i\in I$, $A_i$
denotes the sub-coalgebra indexed by $i$.

It is well known that $A=\underset\longrightarrow {\lim}\,A_i$.
Then $A^*=\underset\longleftarrow {\lim}\,A_i^*$ is a profinite
algebra. If $E$ is a $G$-module (i.e., a linear representation of
$G$) then it is an $A^*$-module. Moreover, if we denote
$E_i=\Hom_{A^*-mod}(A_i^*,E)$, then $E_i$ is an $A_i^*$-module
(acting on $A_i^*$ by the right) and $E=\underset\longrightarrow
{\lim}\, E_i$ as $A^*$-modules. Conversely, if $E$ is an
$A^*$-module such that $E=\underset\longrightarrow {\lim}\, E_i$,
then $E$ is a $G$-module. Moreover, if
$E=\underset\longrightarrow{\lim}\, E_i$ and $ \overline
E=\underset\longrightarrow {\lim}\, \overline E_i$, then
$$\Hom_{G-mod}(E,\overline E)= \Hom_{A^*-mod}(E,\overline E)$$

\begin{definition} Let $E$ be a $G$-module. We shall call
$i$-component of $E$ to
$$E_i=\Hom_{A^*-mod}(A_i^*,E)$$ with the $G$-module
structure induced by the right translations of $G$ on $A_i^*$,
i.e., $g$ acts on $A_i^*$ by $R_{g^{-1}}^{**}$, where $R_g\colon
G\to G$ is the right translation by $g$, $R_g^*\colon A_i\to A_i$
the induced morphism and $R_g^{**}\colon A_i^*\to A_i^*$ the dual
one.
\end{definition}
Note that:
$$E_i=\Hom_{A^*-mod}(A_i^*,E)=\Hom_{G-mod}(A_i^*,E)=(E\underset
k\otimes A_i)^G$$ In particular, the assignation $E\mapsto E_i$
satisfies:\begin{enumerate}
  \item It is functorial, i.e, a morphism of $G$-modules induces a
morphism between its $i$-components.
\item It commutes with base change, i.e.,
$$(E\underset k\otimes
B)_i=E_i\underset k\otimes B$$
 for each base change $k\to
B$.\end{enumerate} Let $E$ be a $G$-module and \[\phi\colon E\to
E\underset k\otimes A=\Hom (G,E)\]
  the structure morphism, i.e., $[\phi (e)](g)=g\cdot e$. This is
a morphism of $G$-modules acting on the latter by the $A$ factor.
By the above said, one has that
\begin{equation}\label{struc} E_i=\phi^{-1}(E\otimes A_i)\end{equation}

\subsection{Classification of principal $G$-bundles}
 Let $G=\Spec A$ be an affine commutative group scheme.
We consider the $G$-module in $A$ given by: $(g\cdot f)(\overline
g)=f( g^{-1}\cdot\overline g)$. Let us denote $G^D$ the dual group
functor of $G$, i.e.,
$$G^D(C)=\Hom_{C-groups}(G_C,(G_m)_C)=\text{ Group of characters of
$G_C$}
$$
  for each $k$-algebra $C$.

 Put as above $A=\underset\longrightarrow {\lim}\, A_i$. Then $\{ A_i^*\}$ is a projective system of
finite commutative algebras and

\begin{proposition}\label{cartier-ind} $G^D=\underset\longrightarrow {\lim}\,\Spec
A_i^*$ (isomorphism of functors).\end{proposition}

\begin{proof} To give an element $\chi_C\in G^D(C)$ is
equivalent to give a character $\chi_C\in A_C$. Since $A_C
=\underset\longrightarrow{\lim}\, A_i\underset k\otimes C $, then
$\chi_C\in A_i\underset k\otimes C$ for some $i$ and $C\cdot\chi$
is a sub-$C$-coalgebra of $A_i\underset k\otimes C$ ; that is,
$\chi_C^*:A_i^*\to C$ is a morphism of $k$-algebras, i.e. an
element of $(\Spec A_i^*)(C)$.
\end{proof} Denoting $Z_i=\Spec A_i^*$, one has then for any
functor $F$
$$\Hom_{func}(G^D, F)=\underset\longleftarrow {\lim}\,
\Hom_{func}(Z_i, F)=\underset\longleftarrow {\lim}\, F(Z_i)$$ For
each $i$, the immersion $Z_i\hookrightarrow G^D$ defines a
character $\chi_i\in A_i\underset k\otimes A_i^*\subset A\underset
k\otimes A_i^*$. Through the isomorphism $ A_i\underset k\otimes
A_i^*=\End_k(A_i^*) $, $\chi_i$ corresponds to the identity of
$A_i^*$.

\begin{definition} The element $\chi_i\in G^D(A_i^*)$ will be
called the universal $i$-character of $G$.\end{definition}

\begin{remarks}\label{remarkcartier}$\,$\newline
\begin{enumerate} \item By Proposition \ref{cartier-ind}  a morphism
of functors $\phi\colon G^D\to F$ is univocally determined by the
images $\phi(\chi_i)$ of the universal $i$-characters of $G$.
\item If $\chi$ is a $C$-valued character, then there exists an
index $i$ such that $\chi$ corresponds to a morphism $f_\chi\colon
\Spec C\to \Spec A_i^*$ and the induced morphism $G^D(A_i^*)\to
G^D(C)$ maps $\chi_i$ onto $\chi$.\end{enumerate}
\end{remarks}

\begin{definition}\label{chi-component} Let $E$ be a $G$-module. For each character $\chi\in G^D(C)$
let $E_{\chi}$ be the sub-$C$-module of $E\underset k\otimes C$
defined as:
$$E_{\chi}=\{\ \overline e\in E\otimes_kC\ :\ g\cdot \overline
e=\chi (g)\overline e\}$$ i.e., $E_{\chi}=(E\underset k\otimes
(C\cdot\chi))^G$ where $C\cdot\chi$ is the sub-$C$-coalgebra of $
A\underset k\otimes C$ generated by $\chi$. We say that $E_\chi$
is the $\chi$-component of $E$.\end{definition}

\begin{example}\label{example} If $E=A$ (ring of functions of $G$), then $A_\chi$
is the $C$-module generated by $\chi^{-1}$: $A_\chi\simeq
C\cdot\chi^{-1}$. Analogously, if $\chi_i$ is the universal
$i$-character, then $(A_i)_{\chi_i}\simeq A_i^*\cdot\chi_i^{-1}$.
\end{example}

\begin{remark}\label{remark} If $\chi\in A_i\underset k\otimes C$, then
$E_{\chi}= (E_i)_{\chi}$. Indeed, from \eqref{struc} one has that
$E_{\chi}\subset E\underset k\otimes A_{\chi}\subset E\underset
k\otimes A_i\underset k\otimes C=(E\underset k\otimes A\underset
k\otimes C)_i$ and then $E_{\chi}= (E_{\chi})_i=(E_i)_{\chi}$.
\end{remark}

\begin{lemma}\label{lemmaB_i} If $\chi_i$ is the universal
$i$-character of $G$, then
$$E_{\chi_i}=\Hom_G(A_i,E)$$ and therefore
$E_{\chi_i}= \Hom_G(A_i,E_i)= \Hom_{A_i^*}(A_i,E_i)$.\end{lemma}

\begin{proof} One has $ E_{\chi_i}=(E_i)_{\chi_i}$ and $(E_i)_{\chi_i}$ is the subspace of
$E_i\underset k\otimes A_i^*=\Hom_k(A_i,E_i)$ defined as $
E_{\chi_i}=\{ f\colon A_i\to E_i, f(g\cdot b)=\chi_i(g)\cdot
f(b)\}$. Now, by definition of $\chi_i$, one has $\chi_i(g)\cdot
e=g\cdot e$ for any $e\in E_i$. Therefore $f\in
E_{\chi_i}\Leftrightarrow f\in
\Hom_G(A_i,E_i)=\Hom_G(A_i,E)$.\end{proof}

\medskip
\noindent{\ \ \it Picard functor}.  Assume now that $Y$ is an
anti-affine scheme with some rational point $p_0$. For each scheme
$Z$ we denote $p_Z\colon Z\to Y\times Z$ the $Z$-valued point
$p_Z(z)=(p_0,z)$. Then the Picard functor of $Y$ is
%
%
$$\mathbf{Pic}(Y)(Z)=\left\{ \aligned &\text{invertible sheaves $\mathcal L$ on
$Y\times Z$}\\ &\text{such that $\mathcal L\res{p_0\times Z}$ is
trivial}\endaligned\right\}$$ Since $Y$ is anti-affine, a morphism
$\lambda\colon\mathcal L\to\mathcal L'$ between invertible sheaves
is univocally determined by the morphism between the fibres at
$p_0$: $\lambda_{p_0}\colon\mathcal L_{p_0}\to\mathcal L'_{p_0}$.
\bigskip

Let $\pi\colon P\to Y$ be a principal $G$-bundle. Since $G$ is
affine, $\pi$ is an affine morphism. Let us denote $\mathcal B =
\pi_*\mathcal O_P$. It is a sheaf of $\mathcal O_Y$-algebras and
$G_Y$-modules. For each character $\chi\in G^D(C)$ let us denote
$\mathcal B_{\chi}$ the $\chi$-component of $\mathcal B$, defined
as in \ref{chi-component}.

\begin{proposition} $\mathcal B_{\chi}$ is an invertible sheaf on $Y_C$.
\end{proposition}

\begin{proof} One has $\mathcal B_{\chi}=(\mathcal B_C\underset
C\otimes (C\cdot \chi))^G$. Hence $\mathcal B_{\chi}$ is stable
under flat base change of $Y$. Then we can assume that $P=G\times
Y$ and then $\mathcal B_{\chi}=\mathcal
O_{Y_C}\cdot\chi^{-1}$.\end{proof} Consequently, a principal
$G$-bundle $\pi\colon P\to Y$ defines a morphism of
 functors of groups:
$$\aligned\phi_{\pi}\colon G^D &\longrightarrow\mathbf{Pic}(Y)\\
\chi &\longmapsto (\pi_*\mathcal O_P)_{\chi}\endaligned$$ and one
has the following:

\begin{theorem}\label{clasFib}{\bf (Classification of principal
$G$-bundles)} Let $Y$ be an anti-affine scheme with some rational
point and $G$ a commutative affine group  scheme. The set
$\operatorname {Prin}(G,Y)$ of isomorphism classes of principal
$G$-bundles over $Y$  is canonically bijective to the set of
morphisms of functors of groups $G^D \to\mathbf{Pic}(Y)$. That is,
the map:
$$\aligned \varphi\colon \operatorname {Prin}(G,Y) &
\longrightarrow \Hom_{groups}(G^D,
\mathbf {Pic}(Y))\\
\pi &\longmapsto \phi_{\pi}\endaligned$$ is bijective.
\end{theorem}

\begin{proof} Let $\phi\colon G^D
\to\mathbf{Pic}(Y)$ be a morphism of functors of groups. One has
to construct, in a functorial way, a sheaf $\mathcal B^{\phi}$ of
$\mathcal O_Y$-$G$-algebras such that $\pi_\phi\colon \Spec
\mathcal B^{\phi}\to Y$ is a principal $G$-bundle. We shall then
see that this construction is the inverse of $\varphi$.

{\sl Construction of $\mathcal B^{\phi}$ as an $\mathcal
O_Y$-$G$-module}: Let $\chi_i$ be the universal $i$-character of
$G$ and let $\mathcal L^{\chi_i}$ be the invertible sheaf on
$Y\times \Spec A_i^*$ (and so a locally free sheaf on $Y$)
corresponding to $\phi (\chi_i)$ and univocally determined by a
fixed isomorphism of $A_i^*$-modules
$$\varphi_i\colon (\mathcal L^{\chi_i})_{p_0}  \overset\sim\to
A_i^*$$ For each inclusion morphism $\Spec
A_i^*\hookrightarrow\Spec A_j^*$ we fix the restriction morphism
$s_{ij}\colon \mathcal L^{\chi_j}\to \mathcal L^{\chi_i}$ as the
only one that coincides with the projection $A_j^*\to A_i^*$ on
the respective fibers over $p_0$. Then one has   $\mathcal
L^{\chi_j}\underset {A^*}\otimes A_i^*= \mathcal L^{\chi_i}$. The
family $\{ \mathcal L^{\chi_i}, s_{ij}\}_i$ is now a projective
system of $O_Y$-modules and $G$-modules. Put $\widehat{\mathcal
L}=\underset\longleftarrow {\lim}\,\mathcal L^{\chi_i}$; one has
$\widehat{\mathcal L}\underset {A^*}\otimes A_i=\mathcal
L^{\chi_i}\underset{A_i^*}\otimes A_i$. Let us denote $$\mathcal
B^{(i)}=\widehat{\mathcal L}\underset{ A^*}\otimes A_i,\qquad
\mathcal B^{\phi}=\underset\longrightarrow {\lim}\,\mathcal
B^{(i)}= \widehat{\mathcal L}\underset{A^*}\otimes A$$ The
isomorphisms $\varphi_i\colon (\mathcal L^{\chi_i})_{p_0}
\overset\sim\to A_i^*$ yield isomorphisms $\mathcal
B^{(i)}_{p_0}\overset\sim\to A_i$ and $\mathcal
B^{\phi}_{p_0}\overset\sim\to A$.

\medskip
{\sl Construction of the algebra structure of $\mathcal B^\phi$}:
\medskip

Let us denote $Z_i=\Spec A_i^*$. For each $i,j$, let $r$ be an
index such that the group structure morphism $m\colon Z_i\times
Z_j\to \underset\longrightarrow {\lim}\, Z_s$ maps into $Z_r$.
Since $\phi$ is a morphism of groups one has:

\begin{equation}\label{algebra} \mathcal L^{\chi_r}\underset
{A_r^*}\otimes (A_i^*\underset k\otimes A_j^*)\simeq\mathcal
L^{\chi_i}\underset k\otimes \mathcal L^{\chi_j}\end{equation}
 and this isomorphism is unique, assuming that, in
the fiber of $p_0$, it coincides with the natural isomorphism
$A_r^*\underset {A_r^*}\otimes (A_i^*\underset k\otimes
A_j^*)=A_i^*\underset k\otimes \mathcal A_j^*$. Now we have a
bilinear morphism:
$$\displaylines{\mathcal B^{(i)}\underset k\otimes\mathcal
B^{(j)}=(\mathcal L^{\chi_i}\underset {A_i^*}\otimes A_i)\underset
k\otimes (\mathcal L^{\chi_j}\underset {A_j^*}\otimes A_j)=
(\mathcal L^{\chi_i}\underset k\otimes \mathcal
L^{\chi_j})\underset {A_i^*\underset k\otimes A_j^* }\otimes
(A_i\underset k\otimes A_j)\overset {\eqref{algebra}}\simeq
\hfill\cr\hfill \mathcal L^{\chi_r}\underset {A_r^*}\otimes
(A_i\underset k\otimes A_j)\longrightarrow \mathcal
L^{\chi_r}\underset {A_r^*}\otimes A_r=\mathcal B^{(r)}}$$
 where
$A_i\underset k\otimes A_j\to A_r$ is the multiplication morphism
on $A$ (which is a morphism of $G$-modules and then of
$A_r^*$-modules). This bilinear morphism is the only  morphism of
$\mathcal O_{Y\times Z_r}$-modules that coincides with the
morphism $A_i\underset k\otimes A_j\to A_r$ at the fibre of
$p_0\times Z_r$. Taking direct limit we have a morphism (of
$G$-modules):
$$\mathcal B^{\phi}\underset{\mathcal O_Y}\otimes\mathcal
B^{\phi}\overset {m^{\phi}}\longrightarrow\mathcal B^{\phi}$$ and
it is the only morphism of $\mathcal O_Y$-$G$-modules that
coincides with the algebra structure morphism $A\underset k\otimes
A\to A$ at the fibre of $p_0$. From the uniqueness of the
construction it is not difficult to see that $m^{\phi}$ gives an
algebra structure on $\mathcal B^{\phi}$ (taking also into account
that it is so for $A\underset k\otimes A\to A$).

Let us denote $P^\phi=\Spec \mathcal B^\phi$. One has a morphism
of $G$-schemes $\pi_\phi\colon P^\phi\to Y$ ($G$ acts trivially on
$Y$). Let us see that $P^\phi\to Y$ is a principal $G$-bundle.
First of all, it is easy to see that the construction of $P^\phi$
is stable under base change. That is, let $f\colon Y'\to Y$ be a
morphism of schemes (and assume that $Y'$ has a rational point
$p'_0$ in the fiber of $p_0$) and let $\phi'\colon G^D\to \mathbf
{Pic}(Y')$ be the morphism of functors obtained by the composition
of $\phi$ with the natural morphism $f^*\colon \mathbf {Pic}(Y)\to
\mathbf {Pic}(Y')$ induced by $f$. Let $\mathcal B^{\phi'}$ the
associated $\mathcal O_{Y'}$-$G$-algebra and $P^{\phi'}=\Spec
\mathcal B^{\phi'}\to Y'$ the associated $G$-scheme over $Y'$.
Then one has a natural isomorphism of $G$-schemes over $Y'$
\[ P^{\phi'}=P^\phi\times_Y Y'\]
Consider now the particular case $Y'=P^\phi$. It is easy to see
that in this case $\phi' (\chi_i)$ is the trivial invertible sheaf
on $Y'\times \Spec A_i^*$. It follows that $\mathcal B^{\phi'}$ is
the trivial $\mathcal O_{Y'}$-$G$-algebra, i.e.,
$P^{\phi'}=Y'\times G$. In other words
\[ P^\phi\times_Y P^\phi = P^\phi\times G\] so $P^\phi\to Y$ is a
principal $G$-bundle.

It remains to prove that the assignations $\pi\mapsto \phi_\pi$
and $\phi\mapsto \pi_\phi$ are inverse to each other.

Let $\phi\colon G^D\to \mathbf {Pic} (Y)$ be a morphism of
functors and $\pi_\phi\colon P^\phi\to Y$ the associated principal
$G$-bundle. Let us see that the morphism of functors associated to
$\pi_\phi$ coincides with $\phi$. By Remark \ref{remarkcartier}
(1), it suffices to see that both  coincide on $\chi_i$. That is,
one has to prove that $\phi(\chi_i)$ is the $\chi_i$-component of
$\mathcal B^\phi$. Recall that $\mathcal B^\phi=\underset\to\lim\,
\mathcal B^{(i)}$, where $\mathcal B^{(i)}=\mathcal
L^{\chi_i}\otimes_{A_i^*}A_i$ and $\mathcal L^{\chi_i}$ is the
invertible sheaf representing $\phi(\chi_i)$. Assume that one has
proved that $\mathcal B^{(i)}$ is the $i$-component of $\mathcal
B^\phi$. Then, by Remark \ref{remark}, $\mathcal B^\phi_{\chi_i} =
\mathcal B^{(i)}_{\chi_i}= (\mathcal
L^{\chi_i}\otimes_{A_i^*}A_i)_{\chi_i} = \mathcal L^{\chi_i}$ (see
Example \ref{example} for the last equality) and we are done. So
let us prove that the $i$-component of $\mathcal B^{\phi}$
coincides with $\mathcal B^{(i)}$. Indeed, locally on $Y$ (for the
Zariski topology), one has $\mathcal L^{\chi_j}\simeq \mathcal
O_Y\underset k \otimes A_j^*$ and then, if $i\leq j$, one has
$\mathcal B^{(j)}\simeq \mathcal O_Y\underset k \otimes A_j$ and
then $(\mathcal B^{(j)})_i=\mathcal B^{(i)}$. Taking direct limit
one concludes.

Now let $\pi\colon P\to Y$ be a principal $G$-bundle and
$\phi_\pi\colon G^D\to \mathbf {Pic} (Y)$ the associated morphism
of functors. We have to prove that $\mathcal B^{\phi_\pi}$ is
canonically isomorphic to $\pi_*\mathcal O_P$ (as $\mathcal
O_Y$-$G$-algebras). Let us denote $\mathcal B=\pi_*\mathcal O_P$.
By definition $\mathcal B^{\phi_\pi} =\underset \to \lim\,
(\mathcal L^{\chi_i}\otimes_{A_i^*}A_i)$, where $\mathcal
L^{\chi_i}$ is the invertible sheaf corresponding to
$\phi_\pi(\chi_i)$, i.e., $\mathcal L^{\chi_i}=\mathcal
B_{\chi_i}$. Since one has a canonical isomorphism of $\mathcal
O_Y$-$G$-modules $\mathcal B_i=\mathcal B_{\chi_i}\underset
{A_i^*}\otimes A_i$ (see Lemma \ref{lema3.7} below) one concludes
that $\mathcal B^{\phi_\pi}$ is canonically isomorphic to
$\mathcal B$ as an $\mathcal O_Y$-$G$-module. From the uniqueness
of the construction of the algebra structure of $\mathcal
B^{\phi_\pi}$ it is not difficult to see that this isomorphism is
in fact an isomorphism of algebras. We are finished.\end{proof}

\begin{lemma}\label{lema3.7} Let
$\pi\colon P\to Y$ be a principal $G$-bundle and $\mathcal
B=\pi_*\mathcal O_P$.  One has a canonical isomorphism of
$\mathcal O_Y$-$G$-modules $$\mathcal B_i=\mathcal
B_{\chi_i}\underset {A_i^*}\otimes A_i$$
%
%
\end{lemma}

\begin{proof} By lemma
\ref{lemmaB_i} one has $\mathcal
B_{\chi_i}=\Hom_{A_i^*}(A_i,\mathcal B_i)$. Hence there is a
natural evaluation morphism:
$$\mathcal B_{\chi_i}\underset {A_i^*}\otimes
A_i=\Hom_{A_i^*}(A_i,\mathcal B_i)\underset {A_i^*}\otimes A_i\to
\mathcal B_i$$ Let us see that it is an isomorphism. After
localizing (for the flat topology) we can assume that $Y=\Spec k$
and $P=G$ and then $\mathcal B=A$ and $\mathcal B_i=A_i$. In this
situation one concludes because
$\Hom_{A_i^*}(A_i,A_i)=\Hom_{A_i^*}(A_i^*,A_i^*)=A_i^*$.
%
%
\end{proof}

\begin{corollary} Under the same hypothesis, every principal
$G$-bundle $P\to Y$ is locally split, i.e., there exists an open
covering $U_i$ of $Y$ such that $P_{\vert U_i}\simeq G\times U_i$.
\end{corollary}

\begin{proof} There exists a ``big enough'' index $j$ such that $G^D$ is
generated by $Z_j$ (as a group). Let $U_i$ be an open covering of
$Y$ trivializing $\mathcal L^{\chi_j}$, i.e., $\mathcal
L^{\chi_j}_{\vert U_i\times Z_j}\simeq \mathcal O_{U_i\times
Z_j}$. Then the composition
$G^D\to\mathbf{Pic}(Y)\to\mathbf{Pic}(U_i)$ is trivial. This
yields that $\mathcal B_{\vert U_i}$ is the trivial $\mathcal
O_{U_i}$-$G$-algebra; that is, $P_{\vert U_i}\simeq U_i\times G$.
\end{proof}

\begin{remark}\label{obse} In the following theorems we shall make use of the following elementary fact:
Let $\chi$ be a $C$-valued character of $G$, i.e., $\chi\in
G^D(C)$. Let $i$ be an index such that $\chi$ corresponds to a
morphism $f_\chi\colon \Spec C\to Z_i$. The induced morphism
$G^D(Z_i)\to G^D(C)$ maps  the universal $i$-character $\chi_i$
onto $\chi$. If $\phi\colon G^D\to \mathbf {Pic} (Y)$ is a
morphism of functors and $\mathcal L^\tau$ denotes the invertible
sheaf representing $\phi(\tau)$ one has
\[ (1\times f_\chi)^*\mathcal L^{\chi_i} =\mathcal L^\chi\] where
$1\times f_\chi\colon Y\times\Spec C\to Y\times Z_i$ is the
morphism induced by $f_\chi$.
\end{remark}

\begin{theorem}\label{antiafin-criterio} Let $G=\Spec A$ be a commutative group and $Y$ an
anti-affine Gorenstein scheme of dimension $g$. Let $\phi\colon
G^D\to \mathbf{Pic}(Y)$ be a morphism and $\pi\colon P\to Y$  the
associated principal $G$-bundle. Put $A=\underset\to\lim\, A_i$,
$\chi_i$ the universal $i$-character, $Z_i=\Spec A_i^*$ and
$\pi_i\colon Y\times Z_i\to Z_i$ the second projection. Then $P$
is anti-affine if and only if
\begin{equation}\label{criterio}   R^g{\pi_i}_*(\omega_Y\underset{\mathcal
O_Y}\otimes\mathcal L^{-\chi_i} )\simeq k_{Z_i}(0)\quad \text{ for
all } i, \end{equation} where $\omega_Y$ is the dualizing sheaf of
$Y$ over $k$, $\mathcal L^\chi$ is the invertible sheaf
representing $\phi(\chi)$  and $k_{Z_i}(0)$ is the ``residual
field of $Z_i$ at the trivial character $0\in G^D(k)$''
 (i.e., $k_{Z_i}(0)=0$ if
$0\not\in  Z_i $ and $k_{Z_i}(0)=k$ if $0\in  Z_i $).
\end{theorem}

\begin{proof} Let us denote $\mathcal
O_{Y\times Z_i}^*=\mathcal Hom_{\mathcal O_Y-mod}(\mathcal
O_{Y\times Z_i},\mathcal O_Y)$. With the same notations as in the
proof of theorem \ref{clasFib}, one has
$$\mathcal B^{(i)}=\mathcal L^{\chi_i}
\underset {A_i^*}\otimes A_i=\mathcal L^{\chi_i}\underset
{\mathcal O_{Y\times Z_i}}\otimes \mathcal O_{Y\times
Z_i}^*=\mathcal Hom_{\mathcal O_{Y\times Z_i}-mod}(\mathcal
L^{-\chi_i},\mathcal O_{Y\times Z_i}^*)$$  Then
$$\mathcal B=\underset\longrightarrow {\lim}\, \mathcal
B^{(i)}=\underset\longrightarrow {\lim}\, \mathcal Hom_{\mathcal
O_{Y\times Z_i}-mod}(\mathcal L^{-\chi_i},\mathcal O_{Y\times
Z_i}^*)$$ and then
$$\displaylines{H^0(P,\mathcal O_P)=H^0(Y,\mathcal
B)=\underset\longrightarrow {\lim}\,H^0(Y,\mathcal B^{(i)}) =\\
\hfill\cr\hfill
  \underset\longrightarrow {\lim}\, H^0(Y\times Z_i, \mathcal Hom_{\mathcal
O_{Y\times Z_i}-mod}(\mathcal L^{-\chi_i},\mathcal O_{Y\times
Z_i}^*))}$$

Since $\mathcal O_{Y\times Z_i}^*$ is the dualizing sheaf of
$Y\times Z_i$ over $Y$, and $\omega_Y\underset{\mathcal
O_Y}\otimes\mathcal O_{Y\times Z_i}^*$ is the dualizing sheaf of
$Y\times Z_i$ over $k$, duality gives
$$ \aligned H^0(Y\times Z_i,\mathcal Hom_{\mathcal O_{Y\times
Z_i}-mod}(\mathcal L^{-\chi_i},\mathcal O_{Y\times Z_i}^*)) &=
H^g(Y\times Z_i,\omega_Y\underset{\mathcal O_Y}\otimes\mathcal
L^{-\chi_i} )^*\\ &= H^0(Z_i,
R^g{\pi_i}_*(\omega_Y\underset{\mathcal O_Y}\otimes\mathcal
L^{-\chi_i} ))^*\endaligned $$
 Hence $P$ is anti-affine if and only if
$$\underset\to\lim\, H^0( Z_i, R^g{\pi_i}_*(\omega_Y\underset{\mathcal O_Y}\otimes\mathcal
L^{-\chi_i} ))^* =k$$ On the other hand, if $i\leq j$, the natural
map $$H^0( Z_i, R^g{\pi_i}_*(\omega_Y\underset{\mathcal
O_Y}\otimes\mathcal L^{-\chi_i} ))^* \to H^0( Z_j,
R^g{\pi_j}_*(\omega_Y\underset{\mathcal O_Y}\otimes\mathcal
L^{-\chi_j} ))^*$$ is injective (use Remark \ref{obse} and
standard properties of the highest direct image). Let $0\in
G^D(k)$ be the trivial character. For any $i$ such that $0\in Z_i$
one has $\mathcal L^{-\chi_i}\otimes_{\mathcal O_{Z_i}}k(0) =
\mathcal L^{-0}=\mathcal O_Y$ (by Remark \ref{obse}). Since the
highest direct image is stable under base change, one obtains that
the fibre of $R^g{\pi_i}_*(\omega_Y\underset{\mathcal
O_Y}\otimes\mathcal L^{-\chi_i} )$ at $0$ is $k$. Moreover one has
a natural epimorphism
$$H^0( Z_i, R^g{\pi_i}_*(\omega_Y\underset{\mathcal
O_Y}\otimes\mathcal L^{-\chi_i} )) \to
R^g{\pi_i}_*(\omega_Y\underset{\mathcal O_Y}\otimes\mathcal
L^{-\chi_i} )\underset{\mathcal O_{Z_i}}\otimes  k(0) = k$$
Putting it all together one concludes.
\end{proof}

\begin{theorem}\label{aaImm} Let $\pi\colon P\to Y$ be a
principal $G$-bundle over  an anti-affine scheme $Y$. If $P$ is
anti-affine then: \begin{enumerate}
\item the associated morphism $\phi \colon G^D\to
\mathbf{Pic} (Y)$ is injective.
\item If $\chi\in G^D(k)$ is a non trivial character, then
$H^0(Y,\mathcal L^{\chi})=0$, where $\mathcal L^{\chi}$ is the
invertible sheaf representing $\phi (\chi)$.\end{enumerate}
\end{theorem}

\begin{proof} (1) If $\chi\in G^D(C)$ is a character in the kernel
of $\phi_\pi$, then $ (\pi_*\mathcal O_P)_{\chi}\simeq \mathcal
O_{Y_C}\cdot \chi $ (as $G_C$-modules). Then $$H^0(P,\mathcal
O_P)\underset k\otimes C=H^0(P_C,\mathcal O_{P_C})\supset C+H^0(Y,
(\pi_*\mathcal O_P)_{\chi})=C+C\cdot\chi.$$  Since $H^0(P,\mathcal
O_P)= k$ , $\chi$ must be trivial.

(2) Let $i$ be an index such that $\chi\in Z_i$. Using Remark
\ref{obse} and Theorem \ref{antiafin-criterio}  one obtains
$$\aligned H^0(Y,\mathcal L^{\chi})=H^g(Y,\omega_Y\otimes\mathcal
L^{-\chi})^* & = H^0(Z_i,R^g{\pi_i}_*(\omega_Y\underset{\mathcal
O_Y} \otimes\mathcal L^{-\chi_i})\underset {\mathcal
O_{Z_i}}\otimes k(\chi))^*\\ & =(k_{Z_i}(0)\underset {\mathcal
O_{Z_i}}\otimes k(\chi))^* =0\endaligned$$
\end{proof}

\medskip
\noindent{\ \ \bf Notations.}  We shall denote
$$\Hom (G^D, \mathbf {Pic}(Y))_0 =\{ \phi\in
\Hom_{groups}(G^D, \mathbf {Pic}(Y)) \text{ satisfying }
\eqref{criterio}\} $$ Then we have proved
$$\operatorname
{Prin}(G,Y)_{ant}=\Hom (G^D, \mathbf {Pic}(Y))_0$$ for any
anti-affine Gorenstein scheme $Y$. If $F,F'$ are two functors of
groups we shall denote by $\Imm_{groups}(F,F')$ the set of
injective morphisms (of functors of groups). We have also proved
that
\[ \operatorname{Prin}(G,Y)_{ant}\subset \Imm_{groups}(G^D,\mathbf{Pic}(Y))\]

\begin{corollary}\label{subAbelian}  An anti-affine principal
$G$-bundle over $Y$ does not admit principal sub-bundles whose
structural group is a strict subgroup $H\subset G$ (strict means
$H\neq G$). In particular, a quasi-abelian variety does not have
strict subgroup schemes with the same abelian part.
\end{corollary}

\begin{proof}
Let $i\colon H\hookrightarrow G$ be a strict subgroup. One has a
surjective and non bijective morphism $i^*\colon G^D\to H^D$. So,
an immersion $G^D\to \bold{Pic}(Y)$ cannot factor through $i^*$.
\end{proof}

Assume now that $Y$ is an abelian variety and denote by
$\operatorname{Prin}(G,Y)_{ant}^{st}$ the set of isomorphism
classes of anti-affine principal $G$-bundles over $Y$ which are
stable under translations on $Y$.

\begin{theorem}\label{clasification-1}
Let $Y$ be an abelian variety, $G$ a connected  commutative affine
group and $\operatorname{Prin}(G,Y)_{ant}^{st}$ the set of
isomorphism classes of anti-affine principal $G$-bundles over $Y$
which are stable under translations on $Y$. Then
\[  \operatorname{Prin}(G,Y)_{ant}^{st}=
\Imm_{groups}(G^D,\mathbf {Pic}^0(Y))
\]
\end{theorem}

\begin{proof} Since $Y$ is an abelian variety, one knows that:
\begin{enumerate}
  \item $\mathbf{Pic}(Y)$ is representable by a smooth scheme,
\item $Y^*=\mathbf {Pic}^0(Y)$ is an abelian variety (the dual
abelian variety of $Y$),
\item if $\mathcal P$ is the Poincar\'e invertible sheaf on $Y\times
{Y^*}$ (the universal one) then $R^g\pi_{Y^*}\mathcal
P=k_{Y^*}(0)$ (and then $R^g\pi_{Y^*}\mathcal P^{-1}=k_{Y^*}(0)$),
  \item $\mathbf {Pic}^I(Y)=\mathbf {Pic}^0(Y)$, where
$\mathbf {Pic}^I(Y)$ is the subgroup-scheme of $\mathbf {Pic}(Y)$
of invertible sheaves that are invariant under translation on $Y$.
  \item $\omega_Y\simeq \mathcal O_Y$.\end{enumerate}

Let $\phi\colon G^D\hookrightarrow \mathbf{Pic}^0(Y)$ be an
injective morphism of functors of groups. Since $\mathbf
{Pic}^I(Y)=\mathbf {Pic}^0(Y)$,  the associated principal
$G$-bundle $\pi\colon P\to Y$ is stable under translations on $Y$.
Moreover $R^g{\pi_i}_*(\omega_Y\underset{\mathcal
O_Y}\otimes\mathcal L^{-\chi_i}) =(R^g\pi_{Y^*}\mathcal
P^{-1})\underset{\mathcal O_{Y^*}}\otimes\mathcal
O_{Z_i}=k_{Y^*}(0)\underset{\mathcal O_{Y^*}}\otimes\mathcal
O_{Z_i}=k_{Z_i}(0)$. By Theorem \ref{antiafin-criterio} $P$ is
anti-affine. Conversely, assume that $P$ is anti-affine. Then
$\phi$ is injective by Theorem \ref{aaImm}. Moreover, if
$\pi\colon P\to Y$ is stable under translations on $Y$, then it is
obvious that each finite sub-scheme $\phi (Z_i)\subset
\mathbf{Pic}(Y)$ is also stable under translations and then
$\phi\colon G^D\to\mathbf{Pic}(Y)$ takes values in $\mathbf
{Pic}^I(Y)=\mathbf{Pic}^0(Y)$.\end{proof}

\subsection{Multiplicative type case}

 Let $G$ be an commutative group of multiplicative type. There
exists a finite Galois extension $K/k$ such that $G_K$ is split
(i.e., it is a diagonalizable $K$-group). Then $G_K^D=X(G_K)$,
i.e., the Cartier-dual functor group is the discrete scheme (over
$K$) associated to the group of characters of $G_K$.  Let us
denote $\mathcal G_{K/k}$ the Galois group of $k\to K$. It is
clear that to give a morphism of functors $G^D\to\bold{Pic}(Y)$ is
equivalent to give a $\mathcal G_{K/k}$-equivariant morphism of
groups $X(G_K)\to\Pic(Y_K)$.

\begin{theorem}\label{multiplicative} If $G$ is a multiplicative type group
and $Y$ is an anti-affine Gorenstein scheme then:
$$\operatorname {Prin}( G, Y)=\Hom_{\mathcal G_{K/k}-groups}(X(G_K),\Pic
(Y_K))$$
 and $$\operatorname {Prin}( G, Y)_{ant}=\Imm_{\mathcal
G_{K/k}-groups}(X(G_K),\Pic_{wd}\, (Y_K))$$ where $\Pic_{wd}\,
(Y_K)=\{$invertible sheaves $\mathcal L$ on $Y_K$ without
associated effective divisors, i.e., such that either $\mathcal
L\simeq\mathcal O_{Y_K}$ or $H^0(Y_K,\mathcal L)=0\}$.
\end{theorem}

\begin{proof} The first equality is due to Theorem \ref{clasFib}
and the isomorphism $G_K^D=X(G_K)$. For the second one, if
$\pi\colon P\to Y$ is an anti-affine principal $G$-bundle, then
the associated morphism $\phi_\pi\colon X(G_K)\to\Pic(Y_K)$ is
injective and takes values in $\Pic_{wd}\, (Y_K)$, by Theorem
\ref{aaImm}. Conversely, if $\phi\colon X(G_K)\to \Pic(Y_K)$ is
injective and takes values in $\Pic_{wd}\, (Y_K)$, then it is easy
to see that the associated principal bundle satisfies conditions
\eqref{criterio} of Theorem \ref{antiafin-criterio} and hence it
is anti-affine.
\end{proof}

\begin{theorem} If $Y$ is an abelian variety and $G$ is a
multiplicative type group then:
$$\operatorname {Prin}(G,Y)_{ant}^{st}= \Imm_{\mathcal
G_{K/k}-groups}(X(G_K),\Pic^0 (Y_K))$$
\end{theorem}

\begin{proof} It follows from Theorem \ref{clasification-1}.
\end{proof}

\subsection{Unipotent case}

 Let $G=\Spec A$ be a commutative affine group scheme and $G_a$ the additive group. We denote
$$\Adit(G)=\Hom_{groups}(G,G_a)$$
 the additive functions
over $G$. It is a vector subspace of $A$.

\begin{proposition} Assume that $\operatorname {char}(k)=p\neq 0$
and let $G$ be a unipotent group with $\dim\, G>0$. Then $\dim_k
\Adit(G)=\infty$.
\end{proposition}

\begin{proof} $\Adit (G_a)=<x, x^p,\dots ,x^{p^n},\dots ,>\subset
k[x]$ is an infinite dimensional vector space. Then, if $\dim\,
G>0$, there exists an epimorphism of groups $f\colon G\to G_a$ and
then $\Adit (G)\supset \Adit (G_a)$, so $\Adit (G)$  has infinite
dimension.
\end{proof}

It is well known that $\Adit (G)$ is canonically isomorphic to the
tangent space $T_e(G^D)$ of $G^D$ at the origin, i.e., the set of
elements of $G^D(k[\varepsilon ])$ that map onto the trivial
element of $G^D(k)$. Moreover, if $U$ is the unipotent part of
$G$, then $\Adit (G)=\Adit (U)$.

\begin{theorem}\label{qaChar-p} Assume $\operatorname {char} (k)>
0$ and let $\pi\colon P\to Y$ be an anti-affine principal
$G$-bundle with $\dim_kH^1(Y,\mathcal O_Y)<\infty$. Then the
unipotent part $U$ of $ G$ is finite. In particular, if $G$ is
quasi-reduced (definition \ref{qreduced}) and connected, then $G$
is a torus.
\end{theorem}

\begin{proof} By Theorem \ref{aaImm}, $\phi_\pi\colon G^D\hookrightarrow
\mathbf{Pic}(Y)$ is injective. Hence
$$T_e(G^D)\to
T_e(\mathbf{Pic}(Y))=H^1(Y,\mathcal O_Y)$$
 is also injective. Then
$\dim_k T_e(G^D)\leq \dim_kH^1(Y,\mathcal O_Y)<\infty$ and $\dim\,
U\leq 0$.

If $G$ is quasi-reduced and connected, then its unipotent part $U$
is finite, quasi-reduced and connected. So $U$ is  a local
rational and finite scheme, i.e., it is trivial. Therefore $G$ is
of multiplicative type  and smooth (because it is quasi-reduced
and connected; see remark \ref{ToroQreduc}).
\end{proof}

If $\operatorname {char}(k)=0$ and $G$ is commutative and
unipotent, then $G\simeq \mathbf E$, where $\mathbf E$ is the
additive group of a finite dimensional vector space $E$, i.e.,
$\mathbf E=\Spec S^\punto_kE^*$.

For any vector space $V$, let us denote $k[V]=S^\punto_kV$ and
$(V)$ the ideal of $k[V]$ generated by $V$. Assume now that
$G\simeq \mathbf E$ and let us denote $A=k[E^*]$ and $A_n=k\oplus
E^*\oplus\cdots\oplus S^n_kE^*$. It is a subcoalgebra of $A$.
Since $\operatorname {char} (k)=0$, using Taylor expansion one can
show that $A_n^*=k[E]/(E)^n$ (isomorphism of algebras) where $
e_1\cdots e_n\in k[E]$ is identified with
$(\frac{\partial}{\partial e_1}\circ\dots\circ
\frac{\partial}{\partial e_n})_0\in A_n^*$. Then:

\begin{proposition} If $\operatorname {char} (k)=0$, then
$$\mathbf E^D=\underset\longrightarrow {\lim}\, \Spec k[E]/(E)^n$$
\end{proposition}

 Let us denote $V=H^1(Y,\mathcal O_Y)$ and
$\mathbf V^*=\Spec S^{\punto }_kV$. Put
$V=\underset\longrightarrow {\lim}\, V_i$, where $V_i$ runs over
the finite dimensional subspaces of $V$. One has
$$\mathbf V^*=\underset\longleftarrow {\lim}\, \mathbf V_i^*$$
 and then
$$(\mathbf V^*)^D=\underset\longrightarrow {\lim}\, (\mathbf
V_i^*)^D$$
 Let $\mathbf{Pic}(Y)_{loc}^0$ be the subfunctor of
groups of $\mathbf{Pic}(Y)$ defined as
$$\mathbf{Pic}(Y)^0_{loc}(C)=\left\{ \aligned & f\colon \Spec
C\to\mathbf{Pic}(Y) \text{ such that $f$ factors trough}
\\ &\text{ some finite, local and rational scheme $\{ Z,z_0\}$:}\\
& \qquad \qquad \xymatrix{ \Spec C \ar[r] ^{f} \ar[rd] _{h}& \mathbf{Pic}(Y)\\
& Z \ar[u]_{g} }\\ &\text{ for some $g\colon Z\to \mathbf{Pic}(Y)$
such that $g(z_0)=0$.} \endaligned \right\} $$ for each
$k$-algebra $C$.

\begin{theorem}\label{PicLoc} Let $V=H^1(Y,\mathcal O_Y)$and $\mathbf V^*=\Spec S^{\punto }_kV$.
One has a canonical isomorphism
$$(\mathbf V^*)^D=\mathbf{Pic}(Y)_{loc}^0$$
\end{theorem}

\begin{proof} By definition of $\mathbf{Pic}(Y)_{loc}^0$ and
taking into account that $(\mathbf V^*)^D=\underset\longrightarrow
{\lim}\, \widetilde Z_i$ with $\widetilde Z_i$ local, rational and
finite schemes, it is enough to show that one has a canonical
isomorphism $\mathbf{Pic}(Y)_{loc}^0(C)=(\mathbf V^*)^D(C)$ for
every local, rational and finite $k$-algebra $C$. Let
   $\mathfrak m\subset C$ be the maximal (nilpotent) ideal. We
have the exact sequence of sheaves of groups on $Y$:
$$0\to \mathfrak m\underset k\otimes \mathcal O_Y\overset
{\text{exp}}\longrightarrow \mathcal O_{Y\times C}^{\text{x}}\to
\mathcal O_Y^{\text{x}} \to 0
$$
 where $B^{\text{x}}$ is the group of invertible elements of
$B$ and $\text{exp} (m\otimes f)=\sum_n\frac 1{n!}\cdot (m\otimes
f)^n$. From the exact sequence of cohomology it follows easily
that:
$$\displaylines{ \mathbf{Pic}(Y)^0_{loc}(C)=H^1(Y,\mathfrak
m\underset k\otimes \mathcal O_Y)=\mathfrak m\underset k\otimes
H^1(Y,\mathcal O_Y)=\underset i{\underset\longrightarrow {\lim}}\,
(\mathfrak m\underset k\otimes V_i)= \hfill \cr \hfill \underset
i{\underset\longrightarrow {\lim}}\, (\underset
n{\underset\longrightarrow {\lim}}\,
\Hom_{k-alg}(k[V_i^*]/(V_i^*)^n,C)) =\underset
i{\underset\longrightarrow {\lim}}\, (\mathbf V_i^*)^D(C)=(\mathbf
V^*)^D(C) \cr}
$$\end{proof}

\begin{theorem}\label{vectorial} Let $Y$ be an anti-affine Gorenstein scheme. If $\operatorname {char} (k)=0$, then
$$\operatorname {Prin}({\mathbf E},Y)=\Hom_{k-lin}(E^*,H^1(Y,\mathcal
O_Y))$$
\end{theorem}

\begin{proof} Denote $V=H^1(Y,\mathcal O_Y)$. By theorems \ref{clasFib} and
\ref{PicLoc} one has
$$\displaylines{\operatorname {Prin}({\mathbf E},Y)=\Hom_{groups}(
\mathbf E^D,\mathbf{Pic}(Y))= \Hom_{groups}( \mathbf
E^D,\mathbf{Pic}_{loc}^0(Y))= \hfill\cr\hfill
\Hom_{groups}(\mathbf
  E^D,(\mathbf V^*)^D)=\Hom_{groups}(\mathbf V^*,
\mathbf E)=\Hom_{k-lin}(E^*,V)}$$
 \end{proof}

Analogously, one has:
\begin{theorem}\label{vectExt} If $Y$ is an abelian variety and $G$
is a
  reduced, connected and commutative unipotent group,
then:

(1) If $\operatorname {char} (k)>0$, then $\operatorname
{Prin}(G,Y)_{ant}^{st}= \operatorname{Quasiabel}(G,Y)= \emptyset$.

(2) If $\operatorname {char} (k)=0$, then $G=\mathbf E$ for some
vector space $E$ and
$$\operatorname {Prin}(G,Y)_{ant}^{st} = \Imm_{k-lin}(E^*,H^1(Y,\mathcal
O_Y))$$
\end{theorem}

\subsection{General case}

Let $G$ be the affine part of a quasi-abelian variety $\cA$. By
Theorem \ref{qaChar-p}, if $\operatorname {char}(k)>0$, then $G$
is a torus. If $\operatorname {char} (k)=0$, then $k$ is a perfect
field and then $G$ is smooth and connected and it splits as a
product $G= U\times\mathcal K$ of its multiplicative type and
unipotent parts. So one has:

\begin{proposition}
If $\cA$ is a quasi-abelian variety, then its affine part
$\cA_{\aff}$ is smooth and it splits as a product $U\times\mathcal
K$, with $U$ a unipotent group and $\mathcal K$ of multiplicative
type.
\end{proposition}

So we assume henceforth that $G$ splits as a product
$G=U\times\mathcal K$, with $U$ a unipotent group and $\mathcal K$
of multiplicative type. Then $G^D=U^D\times \mathcal K^D$. If
$G=\Spec A$ is of multiplicative type, then $A_i^*$ is
geometrically reduced, i.e.,
$$(Z_i)_{\overline k}=\Spec (\overline k\times \overset n
\cdots\times \overline k)$$
 is a discrete finite scheme ($\overline k/k$ being the algebraic
closure). If $G$ is unipotent, then $A_i^*$ is a local $k$-algebra
and then $Z_i$ is a finite and local $k$-scheme. If $G=U\times
\mathcal K$, then $Z_i=Z_i^U\times Z_i^{\mathcal K}=(Z_i)_0\times
(Z_i)_{red}$ where $(Z_i)_0$ is the connected component through
the origin and $(Z_i)_{red}$ is the (geometrically) reduced
sub-scheme of $Z_i$.

\begin{theorem}\label{general} Under the above hypothesis one has: \begin{enumerate}
\item $ \operatorname {Prin}(G,Y) =\operatorname {Prin}(U,Y)\times \operatorname {Prin}({\mathcal
K},Y)$.
\item $\operatorname {Prin}(G,Y)_{ant} =\operatorname {Prin}(U,Y)_{ant} \times \operatorname
{Prin}({\mathcal K},Y)_{ant}$.\end{enumerate} \end{theorem}

\begin{proof} (1) It is immediate because
$$\displaylines{\Hom_{groups}( U^D\times \mathcal K^D,
\mathbf {Pic}(Y))=\hfill\cr\hfill \Hom_{groups}(U^D, \mathbf
{Pic}(Y))\times \Hom_{groups}(\mathcal K^D, \mathbf {Pic}(Y))}$$
(2) We use the anti-affinity criterium of Theorem
\ref{antiafin-criterio}. It is clear that $\mathcal L^{\chi_i}\res
{Z_i^U}= \mathcal L^{\chi_i^U}$ and $\mathcal L^{\chi_i}\res
{Z_i^{\mathcal K}}= \mathcal L^{\chi_i^{\mathcal K}}$. Moreover
$R^g{\pi_i}_*(\omega_Y\underset{\mathcal O_Y}\otimes\mathcal
L^{-\chi_i} )\simeq k_{Z_i}(0)$ if and only if
$R^g{\pi_i}_*(\omega_Y\underset{\mathcal O_Y}\otimes\mathcal
L^{-\chi_i} )\res{(Z_i)_0}\simeq k_{(Z_i)_0}(0)$ and
$R^g{\pi_i}_*(\omega_Y\underset{\mathcal O_Y}\otimes\mathcal
L^{-\chi_i} )\res{(Z_i)_{red}}\simeq k_{(Z_i)_{red}}(0)$. Now,
since the highest cohomology group commutes with base change,
$$R^g{\pi_i}_*(\omega_Y\underset{\mathcal O_Y}\otimes\mathcal
L^{-\chi_i} )\res{(Z_i)_0}=
R^g{\pi_{Z_i^U}}_*(\omega_Y\underset{\mathcal O_Y}\otimes\mathcal
L^{-\chi_i^U} )$$
 and
$$R^g{\pi_i}_*(\omega_Y\underset{\mathcal O_Y}\otimes\mathcal
L^{-\chi_i} )\res{(Z_i)_{red}}=R^g{\pi_{Z_i^{\mathcal
K}}}_*(\omega_Y\underset{\mathcal O_Y}\otimes\mathcal
L^{-\chi_i^{\mathcal K}} )$$
 and we
conclude.\end{proof}

This theorem reduces the computation of principal $G$-bundles (and
anti-affine ones) to the cases when $G$ is either a multiplicative
type or a unipotent group.

\medskip
(*) Let $G$ be a reduced, connected, commutative and affine group
and $\mathcal K$ its multiplicative type part. Let $K/k$ be a
Galois extension such that $(\mathcal K_K)^D$ is discrete. We
denote $\mathcal G_{K/k}$ the Galois group of $k\to K$. Then

\begin{theorem}{\bf (Classification of quasi-abelian
varieties)}\label{classqa} Let $Y$ be an abelian variety, $G$ as
in (*) and $\operatorname{Quasiabel}(G,Y)$ the set of isomorphism
classes of quasi-abelian varieties with affine part isomorphic to
$G$ and abelian part isomorphic to $Y$.
Then
\begin{enumerate}\item If $\operatorname {char} (k)>0$, then
$\operatorname {Quasiabel}(G,Y)\neq \emptyset$ if and only if $G$
is a torus and then:
$$\operatorname {Quasiabel}(G,Y) =\Imm_{\mathcal
G_{K/k}-groups}(X(G_K),\Pic^0(Y_K))/\Aut_{groups}(G\times Y)$$
\item If $\operatorname {char} (k)=0$, then:
$$\displaylines{\operatorname {Quasiabel}(G,Y) =\\
\hfill\cr\hfill \frac{\Imm_{\mathcal
G_{K/k}-groups}(X(G_K),\Pic^0(Y_K))\times
\Imm_{groups}(\Adit(G),H^1(Y,\mathcal O_Y))}{\Aut_{groups}(G\times
Y)}}$$
\end{enumerate}
In another words, to give a quasi-abelian variety $\cA$ with
affine part $G$ and abelian part $Y$ is equivalent to give a
sublattice $\Lambda \subset \Pic^0(Y_K)$, stable under the action
of the Galois group and a linear subspace $V\subset H^1(Y,\mathcal
O_Y)$, up to group automorphisms of $Y$, such that $\Lambda\simeq
X(G_K)$ and $V\simeq \Adit(G)$.
\end{theorem}

A different proof of this result may be found in \cite{Br}. For an
algebraically closed field, this result is given in \cite{Sa01}.

\begin{corollary}\label{corollaryfinite}{\rm (see \cite[Thm. 1]{Ar60} and  \cite[Thm. 4]{Ro61})}
If $k$ is a finite field, then every
quasi-abelian variety is an abelian variety.
\end{corollary}

\begin{proof}
Since $\operatorname {char} (k)>0$ one has that $G_{\aff}$ is a
torus. After base change to $K$ we can assume that it  splits and
then $X(G_{\aff})\simeq \mathbb Z^n$. But $\bold {Pic}^0(Y)$ is a
connected scheme over a finite field, so $\Pic ^0(Y)$ is a finite
set. Therefore $\Imm_{groups} (X(G_{\aff}),\Pic^0(Y))=\emptyset$.
\end{proof}

\end{document}